\documentclass[a4paper,12pt]{article}
\usepackage{graphicx}
\usepackage{amssymb}
\usepackage{amsmath}
\numberwithin{equation}{section}

\begin{document}
\renewcommand{\theequation}{\thesection.\arabic{equation}}
\newcommand{\be}{\begin{equation}}
\newcommand{\ee}{\end{equation}}

\begin{center}

\Large\bf {A Nonlinear Singular Diffusion Equation with Source}
\renewcommand{\thefootnote}{\fnsymbol{footnote}}
\quad\\
\normalsize Pan Jiaqing
\\
Email:jqpan@jmu.edu.cn\\
Department of Mathematics, Jimei University, \\
Xiamen 361021, P.R.China \\
\end{center}

\renewcommand{\thefootnote}{\fnsymbol{footnote}}
\vspace{0.5 cm} \mbox{\bf Abstract:}\quad\rm In this paper, the
existence, uniqueness and dependence on initial value of solution
for a singular diffusion equation with nonlinear boundary
condition are discussed. It is proved that there exists a unique
 global smooth solution which depends on initial data
 continuously.

\mbox{\bf Keywords :}\quad  singular diffusion; global solution;
nonlinear boundary condition

\mbox{\bf 2000 MR Subject Classification: }35K10,35K20,35K60

 \section{Introduction}

~~~~~

In this paper, we consider a boundary value problem
\begin{equation}\left\{
\begin{array}{ll}
u_{t}=(u^{m-1}u_{x})_{x}+u^{p},~~~~~~~~~~~~~~~~~~~~~0<x<1,t>0,\\
u_{x}|_{x=0}=0,~~~~~~~~u_{x}|_{x=1}=-u^{\alpha},~~~~~~~t\geq{0},\\
u|_{t=0}=u_{0},~~~~~~~~~~~~~~~~~~~~~~~~~~~~~~~~~~0\leq{x}\leq{1}.
 \end{array}
\right.
\end{equation}
Where $-1<m<0,~0<p<1,~2-m<\alpha $ and $
0\leq{u_{0}(x)}\leq{M},~\int_{0}^{1}u_{0}(x)dx>0.
$

 The equation in (1.1) arises in many applications in physics
and chemistry . For example, it has been proposed for
$m=\frac{1}{2}$ in plasma physics~(\cite{jgb}), and for $m=-1$ in
the heat conduction in solid hydrogen~(\cite{rose}).

 Although there are many results for $m>0$, the situation is completely different for
$m<0$, where   the equation becomes singular since $u^{m}$ blows
up as $u\longrightarrow{0}$ and
$\int_{0}^{t}\int_{0}^{1}u^{m}dxd\tau$ can be unbounded. Thus
there is essential singularity in (1.1) when $u=0$. Some authors
have discussed the similar problems with $u_{0}>0$. For example,
for positive initial value $u_{0}$,  H. Zhang (\cite{Hz})
discussed the Cauchy problem for $m\in(-1,0]$ with  the conditions
 $$
\lim_{x\longrightarrow{}-\infty}u^{m-1}u_{x}=\lambda,~~~~~~~~~~~~~\lim_{x\longrightarrow{}+\infty}u=1.
 $$
 Where, $\lambda{}>0$.  The author also discussed the first boundary
 value problem for $-1<m<0$ but $u_{0}\geq{}0$ (\cite{pa}).
In order to obtain our conclusions of the paper, we divide the
range $[0,+\infty)$
 into two parts: $[0,t_{\ast}]$ and $[t_{\ast},+\infty)$. We first use Arzela's theorem to
 prove that there exists a function $u^{\ast}$ which solves (1.1) on
 $[0,1]\times[0,t_{\ast}]$. Notice that
 $u^{\ast}(x,t_{\ast})>0$ and  $u^{\ast}(x,t_{\ast})$ is smooth, so we use
 $u^{\ast}(x,t_{\ast})$ as a new initial value and then obtain another solution
 $u^{\ast\ast}$ on  $[0,1]\times[t_{\ast},+\infty)$. Thus we obtain a solution
 $$u(x,t)=\left\{
\begin{array}{ll}
u^{\ast}(x,t),~~~~~~~~~~~~~~~~~~~t\in[0,t_{\ast}],\\
u^{\ast\ast}(x,t),~~~~~~~~~~~~~~~~~~t\in[t_{\ast},+\infty).
\end{array}
\right.
 $$
Finally, with a comparison theorem, we can prove the uniqueness
and the continuous dependence on initial value.

 By a solution of (1.1), we mean a function
 $u(x,t)$ is smooth
enough and satisfies the equation in (1.1), $u_{x}$ is continuous
up to $x=0,1$ and satisfies the boundary condition of (1.1) and
$\lim\limits_{t\longrightarrow{}0}\int_{0}^{1}|u-u_{0}|dx=0.$

The following notations will be used throughout the paper:
$$G_{T}=(0,1)\times(0,T),~~~~~~~~~~~~G=(0,1)\times(0,+\infty),~~~~~~~~~~\overline{u}_{0}=\int_{0}^{1}u_{0}dx.$$

The main results of our paper are as follows:

{\bf Theorem} Assume
\begin{equation}-1<m<0,~0<p<1,~2-m<\alpha,~0\leq{u_{0}(x)}\leq{M},~\overline{u}_{0}>0.\end{equation}
Then there exists a unique global smooth positive solution
$u(x,t)$ to the problem (1.1) such that
$$
u\in{C^{\infty}(G)}\cap{C}([0,+\infty);L^{1}(0,1)).
$$If $u,\hat{u}$ are two solutions
corresponding to $u_{0},\hat{u}_{0}$, then for any  $T>0$, there
is a positive constant $C$ such that
\begin{equation}
\int_{0}^{1}|u-\hat{u}|dx\leq{}C\int_{0}^{1}|u_{0}-\hat{u}_{0}|dx,~~~~~~\mbox{for}~t\in[0,T].
\end{equation}\\

\section{Preliminary lemmas}

~~~~~~~{\bf Lemma 1} Assume $0<u_{0}\leq{}M$ and $u_{0}$ be smooth
enough. For any $T>0$, if $u(x,t)$ is a smooth positive solution
to the problem (1.1) on $G_{T}$, then there exists a positive
constant $C_{0}>0$ such that
$$
\|u\|_{L^{\infty}(G_{T})}\leq{C_{0}},
$$
where, $$C_{0}=[(1-p)T+M^{1-p}]^{\frac{1}{1-p}}.$$

{\bf Proof:} For any $q\geq{0}$, we have
$$
u^{q}u_{t}=u^{q}(u^{m-1}u_{x})_{x}+u^{p+q}.
$$
By Holder's inequality,
\begin{eqnarray}
\frac{1}{1+q}\frac{d}{dt}\int_{0}^{1}u^{q+1}dx&\leq&\int_{0}^{1}u^{p+q}dx\nonumber\\
&\leq&(\int_{0}^{1}u^{1+q}dx)^{\frac{p+q}{1+q}}.
\end{eqnarray}
So,
\begin{eqnarray}
\frac{d}{dt}(\int_{0}^{1}u^{1+q}dx)^{\frac{1-p}{1+q}}&\leq&1-p\nonumber\\
\|u\|_{L^{1+q}(0,1)}&\leq&[(1-p)t+\|u_{0}\|_{L^{1+q}(0,1)}^{1-p}]^{\frac{1}{1-p}}\nonumber\\
&\leq&{C_{0}},~~~~~~~~~~~~~~~~~~\mbox{for}~t\in[0,T],~q\geq{0}.
\end{eqnarray}
By \cite{ad}(Th 2.8, p.25), $\|u\|_{L^{\infty}(G_{T})}\leq{C_{0}}$
.\\

{\bf Lemma 2} Assume $u_{0}$ and $u(x,t)$ be as lemma 1, then
$$
~~~~~~~~~~~~~~~~~~~~~~~~~~~~~~|(u^{\frac{m}{q}})_{x}|\leq{C_{T}}(1+t^{-\frac{1}{2}}),
~~~~~~~~~~~~~~~~~~~~~\mbox{for}~(x,t)\in{[0,1]\times(0,T)}.
$$
Where, $q=\frac{3m-1}{2(m-1)},$ $C_{T}$ depends on $T,m,M,p$
and $\alpha$.\\

{\bf Proof:} Set $ u^{m}=V^{q}, $ then
$$
V_{t}=V^{q-\frac{q}{m}}V_{xx}+(q-1)V^{q-1-\frac{q}{m}}(V_{x})^{2}+\frac{m}{q}V^{1+\frac{qp-q}{m}}.
$$
Differentiating this equation with respect to $x$ and then
multiplying through by $V_{x}$, letting $V_{x}=h$, yields
\begin{eqnarray}
\frac{1}{2}(h^{2})_{t}-V^{q-\frac{q}{m}}hh_{xx}&=&
(3q-2-\frac{q}{m})V^{q-1-\frac{q}{m}}h^{2}h_{x}
+\frac{m}{q}(1+\frac{pq-q}{m})V^{\frac{pq-q}{m}}h^{2}\nonumber\\
&~&+(q-1)(q-1-\frac{q}{m})V^{q-2-\frac{q}{m}}h^{4}.
\end{eqnarray}
For any $0<\tau<T$, let  $\phi(t)$ be a smooth function and
$$\phi(t)=\left\{
\begin{array}{ll}
0,~~~~~~~~~~~~~~~~~~~~~~~~~~~~~~~~~~~t\leq{0},\\
\mbox{monotone},~~~~~~~~~~~~~~~~~~~~~~~0<t<\tau,\\
1,~~~~~~~~~~~~~~~~~~~~~~~~~~~~~~~~~~~t\geq{\tau}.
 \end{array}
\right.
$$
Thus there is a positive constant $C_{\ast}>0$ such that $
0\leq{\frac{d\phi}{dt}}\leq{\frac{C_{\ast}}{\tau}}. $ Set $
Z=(\phi{h})^{2}. $ By \cite{Freid}(Th.6, p.65), we have
$Z\in{}C(\overline{G}_{T})$. Clearly, $ Z|_{t=0}=Z|_{x=0}=0 $ and
(since $\frac{m}{q}-1+\alpha>0$)
\begin{equation}
Z|_{x=1}\leq{(\frac{m}{q})^{2}}C_{0}^{2(\frac{m}{q}-1+\alpha)}.
\end{equation}
Let
$$Z(x_{0},t_{0})=\max\limits_{(x,t)\in{}\overline{G}_{T}}Z(x,t),$$
if $0<x_{0}<1$ and $t_{0}>0$, then
$$
~~~~~~~~~~~~Z_{t}\geq{0},~~~ Z_{x}=0,~~~
Z_{t}-V^{q-\frac{q}{m}}Z_{xx}\geq{0},~~~~~~~~~~\mbox{at}~
(x_{0},t_{0}).
$$
Hence,
$$
~~~~~~~~~~~~~~~~-\phi{\phi_{t}}h^{2}\leq{}\phi^{2}[\frac{1}{2}(h^{2})_{t}-V^{q-\frac{q}{m}}hh_{xx}],~~~~~~~~~
~~\mbox{at}~(x_{0},t_{0}).
$$
Multiplying (2.3) by $\phi^{2}$, we have
$$~~~~~~~~~~~~~
(1-q)(q-1-\frac{q}{m})Z\leq{\frac{m}{q}}(1+\frac{pq-q}{m})u^{p-m+\frac{2m}{q}}\phi^{2}
+\frac{C_{\ast}}{\tau}u^{\frac{2m}{q}+1-m},~~~~~~~~~\mbox{at}~(x_{0},t_{0}).
$$
Since $p<1,m\in(-1,0)$ and $q>0,$ thus
$\frac{m}{q}(1+\frac{pq-q}{m})<0$. Thus we have
$$
~~~~~~~~~~~~~~~~~~~~~~~~~~~~~~~~~~(1-q)(q-1-\frac{q}{m})Z\leq{}\frac{C_{\ast}}{\tau}u^{\frac{2m}{q}+1-m},~~~~~~~~~~~~~~~
~~\mbox{at}~(x_{0},t_{0}) .$$ Notice that$~q=\frac{3m-1}{2(m-1)}$,
hence
$$(1-q)(q-1-\frac{q}{m})>0,~~~~
\frac{2m}{q}+1-m>0.$$
 Let
 $$C^{\ast\ast}=\frac{C_{\ast}C_{0}^{\frac{2m}{q}+1-q}}{(1-q)(q-1-\frac{q}{m})}.$$
 Thus,
\begin{eqnarray}
Z(x_{0},t_{0})&\leq&{}\frac{C_{\ast}C_{0}^{\frac{2m}{q}+1-q}}{\tau(1-q)(q-1-\frac{q}{m})} \nonumber\\
&=&\frac{C^{\ast\ast}}{\tau}.
\end{eqnarray}
Recall from$~Z(x_{0},t_{0})$~that (2.5) holds for all
$(x,t)\in(0,1)\times(0,T)$, specially, for $~0<x<1,~t=\tau~$(here,
$\phi=1,Z =h^{2}(x,\tau)$), thus
$$
~~~~~~~~~~~~~~~~~~~~~~~~~~~~~~~~~~~|h(x,\tau)|\leq{}(\frac{C^{\ast\ast}}{\tau})^{\frac{1}{2}},~~~~~~~~\mbox{for}~(x,\tau)\in{(0,1)\times(0,T)}.
$$
By (2.4), there is another positive constant $C_{T}$ which depends
on $T,m,M,p $ and $\alpha$ such that
\begin{eqnarray}
~~~~~~~~~~~~~~~~~~|h(x,\tau)|&\leq&{}|\frac{m}{q}|C_{0}^{(\frac{m}{q}-1+\alpha)}+(\frac{C^{\ast\ast}}{\tau})^{\frac{1}{2}}\nonumber\\
&\leq&C_{T}(1+\tau^{-\frac{1}{2}}),~~~~~~~~\mbox{for}~(x,\tau)\in{[0,1]\times(0,T)}.\nonumber
\end{eqnarray}
 The proof is complete.\\

We notice that $C_{T}$ increases with respect to $C_{0}$ by (2.5)
and $C_{0}$ increases with respect to $T$ by lemma 1. So we
have\\

 {\bf Corollary } If $T_{1}\leq{T_{2}}$, then $C_{T_{1}}\leq{C_{T_{2}}}$.\\

{\bf Lemma 3} Assume $u_{0}(x)~$and~$u(x,t)~$be  as lemma 1, then
$$
~~~~~~~~~~~~~~~~~~\int_{0}^{1}u(x,t)dx\geq[(\alpha+m-2)t+\int_{0}^{1}u^{2-m-\alpha}_{0}dx]^{\frac{1}{1-m-\alpha}},~
~~~~~\mbox{for}~t\in[0,T].
$$

{\bf Proof:}~Multiplying$~u^{1-m-\alpha}$~to the equation in~(1.1)
yields
$$
\frac{1}{2-m-\alpha}(u^{2-m-\alpha})_{t}=\frac{1}{m}u^{1-m-\alpha}(u^{m})_{xx}+u^{p+1-m-\alpha}.
$$
Because of~$2-m<\alpha$~and $u(x,t)>0$, thus
\begin{eqnarray}
~~~~~~~~~\frac{d}{dt}\int_{0}^{1}u^{2-m-\alpha}dx&=&(2-m-\alpha)[(m-1+\alpha)\int_{0}^{1}u^{-1-\alpha}(u_{x})^{2}dx\nonumber\\
&~&+ \int_{0}^{1}u^{p+1-m-\alpha}dx-1] \nonumber\\
&\leq{}&\alpha+m-2,\nonumber
\end{eqnarray}
so
\begin{equation}
\int_{0}^{1}u^{2-m-\alpha}dx\leq\int_{0}^{1}u_{0}^{2-m-\alpha}dx+(\alpha+m-2)t.
~~~~~~\end{equation} By H\"{o}lder's
inverse-inequality(\cite{ad},~Ch.2, ~Th.2.6), we have
$$
(\int_{0}^{1}udx)^{2-m-\alpha}\leq{\int_{0}^{1}u^{2-m-\alpha}dx}.
$$
Hence by (2.6), we have
\begin{eqnarray}
\int_{0}^{1}u(x,t)dx&\geq&[(\alpha+m-2)t+\int_{0}^{1}u_{0}^{2-m-\alpha}dx]^{\frac{1}{2-m-\alpha}}.\nonumber
\end{eqnarray}
\\

{\bf Lemma 4~} Assume~$u_{1},u_{2}\in{C([0,T],L^{1}(0,1))}$ ~be
two solutions corresponding to ~$u_{10}$ and $u_{20}$, then
$$
~~~~~~~~~~~~~~~~~~~~~\int_{0}^{1}|u_{2}-u_{1}|dx\leq{\int_{0}^{1}|u_{20}-u_{10}|dx}+\int_{0}^{1}\int_{0}^{t}|u_{2}^{p}-u_{1}^{p}|dxd\tau,
~~~~~~~~\mbox{for}~t\in[0,T].
$$

{\bf Proof:~} Take a function $p(x)\in{}C^{\infty}(R)$ such that
$$p(x)=\left\{
\begin{array}{ll}
0,~~~~~~~~~~~~~~~~~~~~~~~~~~~~~~~x\leq{}0,\\
\exp[\frac{-1}{x^{2}}\exp{\frac{-1}{(x-1)^{2}}}],~~~~~~~~~~0<x<1,\\
1,~~~~~~~~~~~~~~~~~~~~~~~~~~~~~~~x\geq{}1.
 \end{array}
\right.
$$
Clearly, $0\leq{}p(x)\leq{}1$ and $p'(x)\geq{}0$. For any given
$\varepsilon>0$, let
$p_{\varepsilon}(x)=p(\frac{x}{\varepsilon}).$ Set
$$
w=\frac{1}{m}(u_{2}^{m}-u_{1}^{m}).
$$
Then $w>0$ iff $u_{2}>u_{1}$. Thus
\begin{eqnarray*}
\int_{0}^{1}(u_{2}-u_{1})_{t}p_{\varepsilon}(w)dx&=&\int_{0}^{1}w_{xx}p_{\varepsilon}(w)dx+\int_{0}^{1}(u_{2}^{p}-u_{1}^{p})p
_{\varepsilon}(w)dx\\
&\leq&(u_{1}^{m-1+\alpha}-u_{2}^{m-1+\alpha})p_{\varepsilon}(w)|_{x=1}+\int_{0}^{1}(u_{2}^{p}-u_{1}^{p})p_{\varepsilon}(w)dx.
\end{eqnarray*}
If $u_{2}(1,t)>u_{1}(1,t)$, then
$(u_{1}^{m-1+\alpha}-u_{2}^{m-1+\alpha})p_{\varepsilon}(w)|_{x=1}<0$
(owing to $\alpha>2-m$). If $u_{2}(1,t)\leq{u_{1}(1,t)}$, then
$w|_{x=1}\leq{}0$ and therefore, $p_{\varepsilon}(w)|_{x=1}=0$.
Thus we always have $
(u_{1}^{m-1+\alpha}-u_{2}^{m-1+\alpha})p_{\varepsilon}(w)|_{x=1}\leq{0}
$ and
\begin{equation}
\int_{0}^{1}(u_{2}-u_{1})_{t}p_{\varepsilon}(w)dx\leq{\int_{0}^{1}(u_{2}^{p}-u_{1}^{p})p_{\varepsilon}(w)dx}.
\end{equation}
Since lemma 3.1 of \cite{est} shows
$$
\int_{0}^{1}(u-\hat{u})_{t}p_{\varepsilon}(w)dx\longrightarrow{}\frac{d}{dt}\int_{0}^{1}[u-\hat{u}]_{+}dx,~~~~~~~~~~~\mbox{as}
~\varepsilon\longrightarrow{0},
$$
thus,
$$~~~~~~~~~~~~~~~\int_{0}^{1}[u_{2}-u_{1}]_{+}dx\leq{\int_{0}^{1}[u_{20}-u_{10}]_{+}dx}+\int_{0}^{1}\int_{0}^{t}[u_{2}^{p}-u_{1}^{p}]_{+}dxd\tau,
~~~~~~\mbox{for}~t\in[0,T],
$$
\begin{equation}\end{equation}
in which, $[u-\hat{u}]_{+}=\max(u-\hat{u},0)$. Similarly,
$$
~~~~~~~~~~~~~\int_{0}^{1}[u_{2}-u_{1}]_{-}dx\leq{\int_{0}^{1}[u_{20}-u_{10}]_{-}dx}+\int_{0}^{1}\int_{0}^{t}[u_{
2}^{p}-u_{1}^{p}]_{-}dxd\tau,~~~~~~\mbox{for}~t\in[0,T],
$$\begin{equation}
\end{equation}
where, $[u-\hat{u}]_{-}=-\min(u-\hat{u},0)$. By~(2.8) and (2.9), we~know that the lemma is true.\\

\section{ Proof of the Theorem}
~~~

We prove our theorem by two steps.\\

 {\bf STEP 1} In this step, we assume that $0<u_{0}\leq{}M$ and $u_{0}$ is smooth enough, $u_{0{}x}|_{x=0}=0,$
 $(u_{0{}x}+u_{0}^{\alpha})|_{x=1}=0$. We will prove that there exists a unique global smooth solution
of
(1.1).\\

For any given $T>0$, we consider the problem (1.1) on
$\overline{G}_{T}$. Make two smooth functions as the following
(\cite{est}, p.997):
$$h(r)=\left\{
\begin{array}{ll}
\frac{1}{2}(2\overline{M})^{m-1},~~~~~~~~~~~~~~~~~~~r\geq{2\overline{M}},\\
\mbox{monotone},~~~~~~~~~~~~~~~~~~~\overline{M}<{r}<2\overline{M},\\
r^{m-1},~~~~~~~~~~~~~~~~~~~~~~~~~~\delta\leq{r}\leq{}\overline{M},\\
\mbox{monotone},~~~~~~~~~~~~~~~~~~~0\leq{r}<\delta,\\
2\delta^{m-1},~~~~~~~~~~~~~~~~~~~~~~~~~r<0.
 \end{array}
\right.
$$
$$g(r)=\left\{
\begin{array}{ll}
\frac{1}{2}(2\overline{M})^{m-2},~~~~~~~~~~~~~~~~~~r\geq{2\overline{M}},\\
\mbox{monotone},~~~~~~~~~~~~~~~~~~\overline{M}<{r}<2\overline{M},\\
r^{m-2},~~~~~~~~~~~~~~~~~~~~~~~~~\delta\leq{r}\leq{}\overline{M},\\
\mbox{monotone}~~~~~~~~~~~~~~~~~~~~0\leq{r}<\delta,\\
2\delta^{m-2}f(r),~~~~~~~~~~~~~~~~~~~r<0.
 \end{array}
\right.
$$
Where,~$0<\delta<\min\limits_{x\in[0,1]}u_{0}(x)$,
$\overline{M}>M$. $\overline{M}$ and $\delta$ are to be
determined. $f(r)\in{}C^{\infty}_{0}(R),~0\leq{}f(r)\leq{}1$ and
$$
f(r)=\left\{
\begin{array}{ll}
1,~~~~~~~~~~~~~~|r|\leq{}1,\\
0,~~~~~~~~~~~~~~|r|\geq{}2.
 \end{array}
\right.
$$

Consider the following problem
 \begin{equation}~~~~~~~~~~~~~~~~~~~~\left\{
\begin{array}{ll}
w_{t}=h(w)w_{xx}+(m-1)g(w)(w_{x})^{2}+w^{p},~~~~~~~~~~0<x<1, t>0,\\
w_{x}|_{x=0}=0,~~~~~~~~~~~w_{x}|_{x=1}=-|w|^{\alpha-1}w,~~~~~~~~~~~t\geq{0},\\
w|_{t=0}=u_{0},~~~~~~~~~~~~~~~~~~~~~~~~~~~~~~~~~~~~~~~~~~~~~~~~~0\leq{x}\leq{1}.
 \end{array}
\right.
\end{equation}

 We first set
$$
\delta=\delta_{0}=\frac{1}{2}\min\limits_{x\in[0,1]}u_{0}(x),
~~~\overline{M}=\overline{M}_{0}=2M.
$$
 The standard parabolic equation theory (\cite{lady},Th.7.4) assumes the existence and uniqueness of
$$w_{0}(x,t)\in{H^{2+\beta,1+\frac{\beta}{2}}}(\overline{G}_{T}),$$
for some $\beta\in(0,1)$, solution of (3.1). By the continuity of
$w_{0}(x,t)$,  there is a $t_{0}>0$, such that
$$
~~~~~~~~~~~~~~~~~~~~~~~~~~~~~~~\delta_{0}\leq{w_{0}}\leq{\overline{M}_{0}},~~~~~~~~~~~~~~~\mbox{for}~(x,t)\in{\overline{G}_{t_{0}}}.
$$
Let
$$
T_{0}=\sup\left\{t_{0}|~~\delta_{0}\leq{}w_{0}\leq{}\overline{M}_{0},~(x,t)\in{}\overline{G}_{t_{0}}\right\}.
$$
Thus by the definition of $h(r)$ and $g(r)$, $w_{0}$ is a solution
of (1.1) on $\overline{G}_{T_{0}}$, or
\begin{equation}~~~~~~~~~~~~~~~~~~~~~~~~~~~~~~~~w_{0}=u,~~~~~~~~~~~~~~~~~~~~~\mbox{for}~t\in[0,T_{0}].\end{equation}
Moreover, $
\lim\limits_{t\longrightarrow{0}}\int_{0}^{1}|u-u_{0}|dx=0. $

Next, we set
\begin{eqnarray*}
\delta=\delta_{1}&=&\frac{1}{2}\min\{[\eta^{\frac{m}{q}}+C_{T}(1+T_{0}^{-\frac{1}{2}})]^{\frac{q}{m}},\delta_{0}\},\\
\overline{M}=\overline{M}_{1}&=&2\max(2M,C_{0}).
\end{eqnarray*}
Where,$$
\eta=[(\alpha+m-2)T+\int_{0}^{1}u^{2-m-\alpha}_{0}dx]^{\frac{1}{1-m-\alpha}}.
$$

For $\delta_{1}$ and $\overline{M}$, there also exists a unique
solution of (3.1)
$w_{1}(x,t)\in{}H^{2+\beta,1+\frac{\beta}{2}}(\overline{G}_{T})$,
and a point $t_{1}$ such that
$$
~~~~~~~~~~~~~~~~~~~~~~~~~~~~~~~\delta_{1}\leq{w_{1}}\leq{\overline{M}_{1}},~~~~~~~~~~~~~~~~~\mbox{for}~(x,t)\in{\overline{G}_{t_{1}}}.
$$
 Let
$$
T_{1}=\sup\left\{t_{1}|~~\delta_{1}\leq{}w_{1}\leq{}\overline{M}_{1},~(x,t)\in{}\overline{G}_{t_{1}}\right\}.
$$
Thus $w_{1}$ is a solution of (1.1) on $\overline{G}_{T_{1}}$, or
\begin{equation}~~~~~~~~~~~~~~~~~~~~~~~~~~~w_{1}=u,~~~~~~~~~~~~~~\mbox{for}~t\in[0,T_{1}].\end{equation}
Clearly, using the lemma 2 of \cite{hon}, we know
$T_{0}\leq{}T_{1}$.

We end this step by showing that~$T_{1}=T$. By the definitions of
$T_{1},~\overline{M}_{1}$ and $\delta_{1}$, there is  a point
$x_{1}\in[0,1]$ such that
\begin{equation}
u(x_{1},T_{1})=\overline{M}_{1},
\end{equation}
or
\begin{equation}
u(x_{1},T_{1})=\delta_{1}.
\end{equation}
 If $T_{1}<T$, then by lemma 1, we have
\begin{eqnarray}
~~~~~~~~~~~~~~~~~~~~~~~~~~~~~~~~~~~~~~~u(x,T_{1})\leq{}C_{0},~~~~~~~~~~~~~~~~\mbox{for}~x\in[0,1].
\end{eqnarray}
Since $ C_{0}<\overline{M}_{1}$, so (3.4)  contradicts (3.6). On
the other hand, since $T_{1}<T$, lemma 3 implies
\begin{eqnarray}
\int_{0}^{1}u(x,T_{1})dx&\geq&[(\alpha+m-2)T_{1}+\int_{0}^{1}u^{2-m-\alpha}_{0}dx]^{\frac{1}{1-m-\alpha}}\nonumber\\
&>&[(\alpha+m-2)T+\int_{0}^{1}u^{2-m-\alpha}_{0}dx]^{\frac{1}{1-m-\alpha}}\nonumber\\
&=&\eta.\nonumber
\end{eqnarray}
Thus there is a $x_{2}\in[0,1]$ such that
$u(x_{2},T_{1})\geq{\eta}$. Using lemma 2 and its Corollary we
have
\begin{eqnarray}
u^{\frac{m}{q}}(x,T_{1})&\leq&{u^{\frac{m}{q}}(x_{2},T_{1})}+C_{T_{1}}(1+{T_{1}}^{-\frac{1}{2}})\nonumber~~~~~~~~~~\\
&\leq&{\eta^{\frac{m}{q}}+C_{T}(1+{T_{0}}^{-\frac{1}{2}})}.\nonumber
\end{eqnarray}
Hence,
\begin{eqnarray}
u(x,T_{1})&\geq&{}[\eta^{\frac{m}{q}}+C_{T}(1+{T_{0}}^{-\frac{1}{2}})]^{\frac{q}{m}}\nonumber~~\\
&\geq&2\delta_{1},~~~~~~~~~~~~~~~~~~~~~~~~~~\mbox{for}~x\in[0,1].
\end{eqnarray}
Clearly, (3.7)  contradicts
 (3.5). Thus, $T_{1}=T$ and
$$
~~~~~~~~~~~~~~~~~~~~~~~~~~~~~~~~~~~~~~w_{1}=u(x,t),~~~~~~~~~~~~~~~~~~~~~~~~~\mbox{for}~(x,t)\in{G_{T}}.
$$
Therefore, $u(x,t)$ is a solution of (1.1) on $G_{T}$. The
 bootstrap argument (\cite{dga}) shows that $u\in{C^{\infty}(G_{T})}$. Recalling from the arbitrariness of
 $T$,
 we know that this step is complete.\\

 {\bf STEP 2} Assume
$u_{0}$ be as (1.2). We will prove that the
conclusions of the theorem are valid.\\

For~$0<\delta<\frac{1}{12}$, let
$$u_{0}^{\ast}=\left\{
\begin{array}{ll}
u_{0},~~~~~~~~~~~~~~~~~~~~~~~~~~x\in[2\delta,1-2\delta],\\
0,~~~~~~~~~~~~~~~~~~~~~~~~~~~x~\overline{\in}~[2\delta,1-2\delta],
 \end{array}
\right.
$$
and
$$
u_{0\delta}=\delta+\delta^{\alpha}x^{2}(1-x)+\int_{0}^{1}u^{\ast}_{0}(y)J(x-\delta{y})dy.
$$
Where, $J$ is a smooth averaging kernel. Clearly, $u_{0\delta}$
satisfies the conditions of STEP 1 and
$$\lim_{\delta\longrightarrow{0}}\|u_{0\delta}-u_{0}\|_{L^{1}(0,1)}=0.$$

For any given $T>0$, we consider the problem
$$~~~~~~~~~~~~~~~~~~~~~~~\left\{
\begin{array}{ll}
u_{t}=(u^{m-1}u_{x})_{x}+u^{p},~~~~~~~~~~~~~~~~~~~~0<x<1, 0<t\leq{}T,\\
u_{x}|_{x=0}=0,~~~~~~u_{x}|_{x=1}=-u^{\alpha},~~~~~~~~0\leq{}t\leq{}T,\\
u|_{t=0}=u_{0\delta},~~~~~~~~~~~~~~~~~~~~~~~~~~~~~~~~~0\leq{x}\leq{1}.
 \end{array}
\right. $$ STEP 1 assures that there is a smooth
solution~$u_{\delta}\in{C^{\infty}}(G_{T})\bigcap{C}([0,T];
L^{1}(0,1))$ and
$$
~~~~~~~~~~~~~~\int_{0}^{1}u_{\delta}(x,t)dx\geq{\int_{0}^{1}u_{0\delta}(x)dx}-\int_{0}^{t}u^{m-1+\alpha}_{\delta}(1,\tau)d\tau,
~~~~~~~~~~~~~~~\mbox{for}~t\in[0,T].$$
\begin{equation}\end{equation}
Recalling from
$\int_{0}^{1}u_{0\delta}dx\longrightarrow\overline{u}_{0}$ as
$\delta\longrightarrow{}0$, hence we know that there are
$\delta_{0}$ and $t_{0}$ such that
\begin{equation}
~~~~~~~~~~~~~~~~~~~~~~~~~~~~~\int_{0}^{1}u_{\delta}(x,t)dx\geq{\frac{1}{2}}\overline{u}_{0},~~~
~~~~~~~~~\mbox{for}~\delta\in(0,\delta_{0}),~t\in[0,2t_{0}].
\end{equation}
For any given $\tau\in(0,2t_{0}]$, lemma 1 and lemma 2 and
~Arzela's theorem assure the existence of subsequence
$\left\{u_{\delta_{k}}(x,t)\right\}$ and a function
$u^{\ast}(x,t)$ such that
\begin{equation}~~~~~~~~~~~~~~~~~~~~~~~~\lim_{\delta_{k}\longrightarrow{0}}u_{\delta_{k}}(x,t)=u^{\ast}(x,t),~~~
~~~~~~\mbox{uniformly on}~x\in[0,1]\end{equation} for
$t\in[\tau,2t_{0}]$. On the other hand, (3.9) implies that for
any~$\delta\in(0,\delta_{0})$, there is a point$~(x_{3},t)$ such
that
\begin{equation}
~~~~~~~~~~~~~~~~~
u_{\delta}(x_{3},t)\geq{\frac{1}{2}}\overline{u}_{0},~~~~~~~~~~\mbox{for}~t\in[\tau,2t_{0}].
\end{equation}
By lemma 2,
\begin{eqnarray*}~~~~~~~~~~~~~~~~~~u_{\delta}^{\frac{m}{q}}(x,t)&\leq&{u^{\frac{m}{q}}_{\delta}}(x_{3},t)+C_{T}(
1+t^{-\frac{1}{2}})\\
&\leq& {u^{\frac{m}{q}}_{\delta}}(x_{3},t)+C_{T}(
1+\tau^{-\frac{1}{2}}),~~~~~~~~\mbox{for}~(x,t)\in[0,1]\times[\tau,2t_{0}].
\end{eqnarray*}
Using (3.11),
\begin{eqnarray}
~~~~u_{\delta}(x,t)&\geq{}&[(\frac{\overline{u}_{0}}{2})^{\frac{m}{q}}+C_{T}(1+\tau^{-\frac{1}{2}})]
^{\frac{q}{m}}\nonumber\\
&>&0,~~~~~~~~~~~~~~~~~~~~~~~\mbox{for}~(x,t)\in[0,1]\times[\tau,2t_{0}].
\end{eqnarray}
Set
\begin{eqnarray*}
A&=&u_{\delta}^{m-1},\\
B&=&\frac{4(m-1)}{m^{2}}((u_{\delta}^{\frac{m}{2}})_{x})^{2}+u_{\delta}^{p}.
\end{eqnarray*}
Thus lemma 2 and (3.12) imply that there is a~positive constant
$\mu$ which doesn't depend on $\delta\in(0,\delta_{0})$ such that
~
$$
~~~~~~~~~~~~~~~~~~~~~~~~0<A<\mu,~~~~~ |B|<\mu,
~~~~~~~~~~~\mbox{for}~(x,t)\in[0,1]\times[\tau,2t_{0}].
$$
Notice that $u_{\delta}$~satisfies the linear equation
$$
\frac{\partial}{\partial{t}}u_{\delta}=A\frac{\partial^{2}}{\partial{x^{2}}}u_{\delta}+
B.
$$
For any $\varepsilon\in(0,\frac{1}{2})$, \cite{BHG}(p.104)~shows
that there are positive constants $h$,$\nu$~and~$C$, which don't
depend on $\delta\in(0,\delta_{0})$,  such that
$$
|u_{\delta}(x,t_{2})-u_{\delta}(x,t_{1})|\leq{C}|t_{2}-t_{1}|^{h},
$$
for$~t_{1},~t_{2}\in[\tau, 2t_{0}],~|t_{1}-t_{2}|<\nu,
x\in[\varepsilon,1-\varepsilon]$. Certainly, we also have
$
|u_{\delta_{k}}(x,t_{2})-u_{\delta_{k}}(x,t_{1})|\leq{C}|t_{2}-t_{1}|^{h}.
$
Letting~$\delta_{k}\longrightarrow{0}$ yields
$$
|u^{\ast}(x,t_{2})-u^{\ast}(x,t_{1})|\leq{C}|t_{2}-t_{1}|^{h},
$$
for$~t_{1},~t_{2}\in[\tau, 2t_{0}],~|t_{1}-t_{2}|<\nu,
x\in[\varepsilon,1-\varepsilon]$. Thus, for any given $x\in(0,1),$
$u^{\ast}(x,t)$~is~continuous with respect to $t\in[\tau,2t_{0}]$.
On the other hand, lemma 2 implies that there is a positive
constant $K$ such that $|u_{\delta_{k}}|\leq{}K$ on
$(x,t)\in[\tau,2t_{0}]\times[0,1]$, so
$|u_{\delta_{k}}(x_{2},t)-u_{\delta_{k}}(x_{1},t)|\leq{}K|x_{2}-x_{1}|$.
Letting $\delta_{k}\longrightarrow{}0$, we have
$$
~~~~~~~~~~~~~~~~~|u^{\ast}(x_{2},t)-u^{\ast}(x_{1},t)|\leq{}K|x_{2}-x_{1}|~~~~~~~\mbox{uniformly
on}~x_{1},x_{2}\in[0,1],t\in[\tau,2t_{0}].
$$
Now we have~$u^{\ast}\in{C([0,1]\times[\tau,2t_{0}])}$ and
~$u^{\ast}(x,t)>0$ for $(x,t)\in{}[0,1]\times[\tau,2t_{0}]$. By
lemma 5 of \cite{Gildi}, we know that $u^{\ast}~$ satisfies the
equation and the boundary conditions of (1.1). Clearly,
$u^{\ast}\in {C([\tau,2t_{0}];~ L^{1}(0,1))}$. Because $\tau>0$ is
arbitrary, so  $u^{\ast}\in {C((0,2t_{0}]; L^{1}(0,1))}$.

To show that $u^{\ast}$~is a solution of~(1.1)~on~$G_{2t_{0}}$, we
want to prove $\|u^{\ast}-u_{0}\|_{L^{1}(0,1)}\longrightarrow{}0$
as $t\longrightarrow{}0$.

For any ~$\delta_{k}, \delta_{k+j}$, lemma 4 implies
\begin{eqnarray*}
\|u_{\delta_{k}}-u_{\delta_{k+j}}\|_{L^{1}(0,1)}&\leq&\|u_{0\delta_{k}}-u_{0\delta_{k+j}}\|_{L^{1}(0,1)}
+\int_{0}^{t}\|u^{p}_{\delta_{k}}-u^{p}_{\delta_{k+j}}\|_{L^{1}(0,1)}d\tau\\
&\leq&\|u_{0\delta_{k}}-u_{0}\|_{L^{1}(0,1)}+\|u_{0_{\delta_{k}+j}}-u_{0}\|_{L^{1}(0,1)}\\
&~&
+\int_{0}^{t}\|u^{p}_{\delta_{k}}-u^{p}_{\delta_{k+j}}\|_{L^{1}(0,1)}d\tau,~~~~~~~~~\mbox{for}~
t\in(0,2t_{0}].
\end{eqnarray*}
Letting  $j\longrightarrow{\infty}$ yields
$$
~~~~~~~~~~~~~~~~~~~\|u_{\delta_{k}}-u^{\ast}\|_{L^{1}(0,1)}\leq\|u_{0\delta_{k}}-u_{0}\|_{L^{1}(0,1)}
+\int_{0}^{t}\|u^{p}_{\delta_{k}}-u^{\ast p
}\|_{L^{1}(0,1)}d\tau,~~~~~~~\mbox{for}~ t\in(0,2t_{0}].
$$
Notice that
\begin{eqnarray*}
\|u^{\ast}-u_{0}\|_{L^{1}(0,1)}&\leq&\|u^{\ast}-u_{\delta_{k}}\|_{L^{1}(0,1)}
+\|u_{\delta_{k}}-u_{0\delta_{k}}\|_{L^{1}(0,1)}+\|u_{0\delta_{k}}-u_{0}\|_{L^{1}(0,1)}\\
&\leq&2\|u_{0\delta_{k}}-u_{0}\|_{L^{1}(0,1)}+\|u_{{\delta_{k}}}-u_{0\delta_{k}}\|_{L^{1}(0,1)}\\
&~& +\int_{0}^{t}\|u^{p}_{\delta_{k}}-u^{\ast
{p}}\|_{L^{1}(0,1)}d\tau,~~~~~~~~~~~~~~~~~~~~~~\mbox{for}~
t\in(0,2t_{0}].
\end{eqnarray*}
Thus,
$$
\lim_{t\longrightarrow{0}}\|u^{\ast}-u_{0}\|_{L^{1}(0,1)}\leq{}2\|u_{0\delta_{k}}-u_{0}\|_{L^{1}(0,1)}.
$$
Letting~$\delta_{k}\longrightarrow{0}$ shows
$$
\lim_{t\longrightarrow{0}}\|u^{\ast}-u_{0}\|_{L^{1}(0,1)}=0.
$$

 Next, we consider the problem
\begin{equation}\left\{
\begin{array}{ll}
u_{t}=(u^{m-1}u_{x})_{x}+u^{p},~~~~~~~~~~~~~~~~~~~~~~~0<x<1,~~t_{0}<t\leq{T},\\
u_{x}|_{x=0}=0,~~~~~~~~~u_{x}|_{x=1}=-u^{\alpha},~~~~~~~~t_{0}\leq{}t\leq{}T,\\
u|_{t=t_{0}}=u^{\ast}(x,t_{0}),~~~~~~~~~~~~~~~~~~~~~~~~~~~~0\leq{x}\leq{1}.
 \end{array}
\right.
\end{equation}
Since~$u^{\ast}(x,t_{0})>0$ and $u^{\ast}(x,t_{0})$ is smooth
enough for $x\in[0,1]$, the conclusion of STEP 1 shows that there
is a function $u^{\ast\ast}$ to solve (3.13). Now we define a
function
$$u(x,t)=\left\{
\begin{array}{ll}
u^{\ast},~~~~~~~~~~~~~~~~~~~~~~~~~~~~~ t\in[0,t_{0}],\\
u^{\ast\ast},~~~~~~~~~~~~~~~~~~~~~~~~~~~~t\in[t_{0},T].
 \end{array}
\right.
$$
Clearly, $u~$ is a solution of~(1.1)~in~$G_{T}$ and~the
bootstrap~argument (\cite{dga}) shows $u\in{C^{\infty}}(G_{T})$.

 To end the proof of our theorem, we assume
$$\overline{u}(x,t)=\left\{
\begin{array}{ll}
u_{11},~~~~~~~~~~~~~~~~~~~~~~~~~~~~~~~ t\in[0,t_{\ast}],\\
u_{12},~~~~~~~~~~~~~~~~~~~~~~~~~~~~~~~t\in[t_{\ast},T],
 \end{array}
\right.
$$
$$\overline{\overline{u}}(x,t)=\left\{
\begin{array}{ll}
u_{21},~~~~~~~~~~~~~~~~~~~~~~~~~~~~~~~ t\in[0,t_{\ast}],\\
u_{22},~~~~~~~~~~~~~~~~~~~~~~~~~~~~~~~t\in[t_{\ast},T],
 \end{array}
\right.
$$
in which, $\overline{u}$ and $\overline{\overline{u}}$ are two
solutions corresponding to initial values $u_{10}$ and $u_{20}$.
Thus lemma 4 shows
$$
\|u_{21}-u_{11}\|_{L^{1}(0,1)}\leq\|u_{20}-u_{10}\|_{L^{1}(0,1)}+\int_{0}^{t}\|u^{p}_{21}-u^{p}_{11}\|
_{L^{1}(0,1)}d\tau,~~~~~~~~\mbox{for}~t\in[0,t_{\ast}].
$$
By lemma 1,~$u_{ij}$ are bounded on $G_{T}$ for ~$i,j=1,2$. Hence
we can set ~$t_{\ast}$ small enough such that
\begin{equation}
~~~~~~~~~~~~~~~~~~~~~~~~~\|u_{21}-u_{11}\|_{L^{1}(0,1)}\leq2\|u_{20}-u_{10}\|_{L^{1}(0,1)},~~~~~~~~~~\mbox{for}~t\in[0,t_{\ast}].
\end{equation}
Notice that (2.7) yields
\begin{eqnarray*}
~~~~~~~~~~~~~~~~~~~\frac{d}{dt}\int_{0}^{1}|u_{22}-u_{12}|dx&\leq&{\int_{0}^{1}}|u^{p}_{22}-u^{p}_{12}|dx\\
&\leq&p\xi^{p-1}\int_{0}^{1}|u_{22}-u_{12}|_{L^{1}(0,1)}dx,~~~~~~~~~~~\mbox{for}~t\in[t_{\ast},T],
\end{eqnarray*}
in which,
\begin{eqnarray*}
~~~~~~~~~~~~~\xi&=&{\min_{(x,t)\in[0,1]\times[t_{\ast},T]}}(u_{12},u_{22})\\&>&0.
\end{eqnarray*}
 Using (3.14) we have
\begin{eqnarray*}
~~~~~~~~~~~~~~~~\|u_{22}-u_{12}\|_{L^{1}(0,1)}&\leq&(\|u_{22}-u_{12}\|_{L^{1}(0,1)})_{t=t_{\ast}}e^{p\xi^{p-1}t}\\
&\leq&2\|u_{20}-u_{10}\|_{L^{1}(0,1)}e^{p\xi^{p-1}t},~~~~~~~~~~~~~~~~~~\mbox{for}~t\in[t_{\ast},T].
\end{eqnarray*}\begin{equation}\end{equation}
It follows from $0<p<1$ that $e^{p\xi^{p-1}t}\leq1$. Combining
(3.14) and (3.15) yields (1.3), the uniqueness of the solution is
followed immediately.
\\

The author is pleased to express his gratitude to  Prof. Li
Ta-tsien for his valuable guidance.

\end{document}